\documentclass[12pt]{amsart}

\usepackage{amsmath, amsfonts, amssymb}

\newtheorem{thm}{Theorem}[section]
\newtheorem{lem}[thm]{Lemma}
\newtheorem{prop}[thm]{Proposition}

\usepackage{tikz}

\usetikzlibrary{decorations.markings}
\usetikzlibrary{shapes.geometric}

\pgfdeclarelayer{edgelayer}
\pgfdeclarelayer{nodelayer}
\pgfsetlayers{edgelayer,nodelayer,main}

\tikzstyle{none}=[inner sep=0pt]
\definecolor{hexcolor0xf81e1c}{rgb}{0.973,0.118,0.110}
\definecolor{hexcolor0x3c00ff}{rgb}{0.235,0.000,1.000}

\tikzstyle{blackvertex}=[circle,fill=black,draw=black, scale=0.61]
\tikzstyle{vertex}=[circle, fill=none,draw=black, scale=0.61]

\tikzstyle{box}=[rectangle,fill=none,draw=none]

\tikzstyle{arc}=[-,draw=black,postaction={decorate},decoration={markings,mark=at position .8 with {\arrow{>}}},line width=0.61]

\begin{document}
\title{Obstructions to some injective oriented colourings}
\author{Russell J.~Campbell \and Nancy E.~Clarke
\and Gary MacGillivray}

\begin{abstract}
Each of several possible definitions of local injectivity for a homomorphism of an oriented graph $G$ to an oriented graph $H$ leads to an injective oriented colouring problem.
For each case in which such a problem is solvable in polynomial time, we identify a 
set $\mathcal{F}$ of oriented graphs such that an
oriented graph $G$ has an injective oriented colouring with the given 
number of colours if and only if there is no $F \in \mathcal{F}$
for which there is a locally-injective homomorphism of $F$ to $G$.
\end{abstract}

\maketitle

\section{Introduction}
\label{IntroSec}
\subsection{Motivation}
A $k$-colouring of a graph $G$ is called \emph{locally-injective} if, for every $x$ in $V(G)$,  no two vertices in $N(x)$ are assigned the same colour.  
The definition allows for adjacent vertices to be assigned the same colour.
Thus these may not be proper colourings.
A proper locally-injective colouring of a graph $G$ is the same as a proper colouring of $G^2$,
the graph with the same vertex set as $G$, and $xy \in E(G^2)$ whenever the distance in $G$ between $x$ and $y$ is at most two.

Locally-injective colourings of undirected graphs were first explicitly studied by Hahn, Kratochvil, Si\v{r}an and Sotteau \cite{gena}.  
Subsequent papers have considered various graph classes including 
chordal graphs \cite{hrs08}, and planar graphs (see \cite{KS}).
Also see \cite{cky12, dhr10}.
Locally-injective list colourings are studied in \cite{bi12, fk06, LX12}.  
The closely related topic of the complexity of locally-injective homomorphisms of graphs has been extensively investigated by Fiala, Kratochvil, 
and others (e.g. see \cite{blt12,fkp08,fp10,fpt08,lt11, R14}).  
Other related topics include locally-surjective homomorphisms and locally-bijective homomorphisms (see \cite{bkm11,fpt08}).

We seek to continue the study of locally-injective colourings of oriented graphs, which
are defined and discussed below.
These
first arose as an example in the work of Courcelle \cite{courcelle}, where they were called \emph{semi-strong} colourings
(also see \cite{RS94}).
Our main goal is the following:
in cases where the locally-injective colouring problem is solvable in polynomial time,
we seek to identify the obstructions, or forbidden configurations, which prevent the
given oriented graph from being coloured with the pallette of colours being used.

\subsection{Homomorphisms and oriented colourings}

A homomorphism of a directed graph $G$ to a directed graph $H$ is a 
function $f: V(G) \to V(H)$ such that if $xy \in E(G)$, then $f(x)f(y) \in E(H)$.
The definition is the same if $G$ and $H$ are graphs.
We write 
$G \to H$ to denote the existence of a homomorphism of $G$ to $H$.
The book \cite{hn_book} contains a wealth of information about homomorphism of 
graphs and digraphs.

A graph or directed graph is called \emph{reflexive} if there is a loop at every vertex, and 
is called \emph{irreflexive} if no vertex has a loop.

It follows on comparing the definitions that a proper $k$-colouring of a graph $G$ is a homomorphism
of $G$ to the graph $H$ with vertex set equal to the set of colours and $c_1c_2 \in E(H)$ if there exists
$xy \in E(G)$ such that $f(x) = c_1$ and $f(y) = c_2$.  
Thus if $G$ has a proper $k$-colouring, then there exists an irreflexive graph $H$ with $k$ vertices such that $G \to H$.
Conversely, if there exists an irreflexive graph $H$ with $k$ vertices such that $G \to H$, then
$G$ has a proper $k$-colouring: adjacent vertices of $G$ are assigned adjacent (hence different) images in $H$.
Since there is a homomorphism of any such graph $H$ to the complete graph $K_k$, and a composition of homomorphisms is a homomorphism, 
it follows that a graph $G$ has a proper $k$-colouring if and only if $G \to K_k$. 

A directed graph $G$ is an \emph{oriented graph} if, for any two different vertices 
$x$ and $y$, at most one of $xy, yx$ is an arc of $G$.  
The oriented graphs we consider may have loops.
An oriented graph $G$ can be viewed as arising from a  graph $G'$ without multiple edges 
by assigning a direction, or \emph{orientation}, to each edge.
The graph $G'$ is called the \emph{underlying graph} of $G$, and $G$ is referred to as \emph{an orientation} of $G'$.

Let $G$ be an oriented graph.
A function $f: V \to \{1, 2, \ldots, k\}$ is an \emph{oriented $k$-colouring of $G$} if,
 for all $xy, vw \in E$, if $f(x) = f(w)$, then $f(y) \neq f(v)$.
The condition to be satisfied can be regarded as requiring the colour assignment to respect the orientation of the edges:  
if there is an arc from a vertex assigned colour $i$ to a vertex assigned colour $j$, then there is no arc
from a vertex of colour $j$ to a vertex of colour $i$;
either every arc between $f^{-1}(i)$ and $f^{-1}(j)$ is oriented from the vertex assigned colour $i$ to the vertex 
assigned colour $j$, 
or every arc between $f^{-1}(i)$ and $f^{-1}(j)$ is oriented from the vertex assigned colour $j$ to the vertex 
assigned colour $i$.
An oriented colouring, $f$, is called \emph{proper} if $f(x) \neq f(y)$ whenever $x \neq y$ and $xy \in E$.
Oriented colouring is a well-studied topic; 
see \cite{survey}.

As above, it follows on comparing the definitions that an irreflexive oriented graph $G$ has a 
proper oriented $k$-colouring if and only if there exists an irreflexive oriented graph $H$ with $k$ 
vertices such that $G \to H$.
Since any such oriented graph $H$ has a homomorphism to a tournament on $k$ vertices, 
if follows that an irreflexive oriented graph $G$ has a proper oriented $k$-colouring if and only if
it has a homomorphism to some tournament on $k$ vertices.
Similarly, an oriented $k$-colouring (not necessarily proper) of an oriented graph $G$ is the same as a homomorphism
to a reflexive tournament on $k$ vertices (i.e.~to a tournament with a loop at each vertex).

\subsection{Duality theorems}
Let $H$ be a (directed) graph.
By a \emph{duality theorem} for homomorphism to $H$ we mean a theorem of the following form:
a (directed) graph $G$ has a homomorphism to $H$ if and only if no element of a particular set $\mathcal{F}$
of (directed) graphs has a homomorphism to $G$.
The elements of $\mathcal{F}$ are called \emph{obstructions} to homomorphism to $H$.
An example of such a theorem is that a directed graph $G$ has a homomorphism to $T_n$, the transitive tournament
on $n$ vertices, if and only if the directed path $P_{n+1}$ does not have a homomorphism to $G$.  
It follows that $G \to T_n$ if and only if $G$ contains neither a directed cycle nor a directed path on $n+1$ vertices.

When $H$ is finite there is a succinct certificate that $G \to H$, namely a description of the homomorphism.
Duality theorems often make it possible to give a succinct certificate that there is
no homomorphism of $G$ to $H$, namely a description of homomorphism of an obstruction to $G$;
 for example, this is possible when there is guaranteed to exist an obstruction of size polynomial in $|V(G)|$.
Decision problems for which there is a succinct certificate for either answer 
are in NP $\cap$ co-NP, and such problems tend to be in P.
If it is not difficult to test whether an element of the set $\mathcal{F}$ has
a homomorphism to  $G$, then a duality theorem leads to a polynomial time algorithm for
testing whether $G \to H$.  

A directed graph $H$ has \emph{bounded treewidth duality} if there exists a positive integer $w$
such that every element of the set $\mathcal{F}$ has treewidth at most $w$.
Hell, Ne\v{s}et\v{r}il and Zhu \cite{HNZ} have shown that if $H$ has bounded treewidth 
duality, then a consistency check algorithm (see Section \ref{Sec4} for an example) can be used to determine 
whether there is a homomorphism for a given directed graph $G$ to $H$.
In such situations, duality theorems are sometimes found by analyzing the situations where a consistency
check algorithm fails and tracing backwards; this is the approach we take to find Theorem \ref{ios-refl-2}.

There are duality theorems for colouring problems.
The statement above in the first paragraph of Section 1.3 with $n = 2$ is a duality theorem for proper oriented 2-colouring.
A dichotomy theorem for proper oriented 3-colouring is as follows.
The \emph{net-length} of an oriented cycle equals the number of forwards arcs minus the number of backwards
arcs.  
Let $\mathcal{F} = \{P_4 \cup C: C$ is an oriented cycle of net-length not divisible by 3$\}$. 
Then, 
an oriented graph $G$ has a proper oriented $3$-colouring if and only if no element of $\mathcal{F}$ 
has a homomorphism to $G$.
This statement is derived by combining the duality theorems for homomorphism to each of the two
tournaments on three vertices.

For more information on homomorphisms and duality, see \cite{hn_book}.

\subsection{Injective oriented colourings}

When some vertex of the directed graph $H$ has a loop, 
any directed graph $G$ has a homomorphism to $H$: 
map all vertices of $G$ to a vertex of $H$ with a loop.  
Thus, when loops are present, the existence of a homomorphism 
of $G$ to $H$ is a non-trivial question only in the presence of some side condition like selecting the image of each vertex from a list of possible images, or local injectivity.  

Three natural definitions of local-injectivity of a homomorphism $f$ from an \emph{input} directed graph $G$ to a \emph{target} directed graph $H$ are as follows:
for every vertex $x \in V(G)$, the function $f$ is injective when restricted to

\vspace{-2.5mm}
\begin{enumerate}
\item the  \emph{in-neighbourhood}, $N^-(x) $ (equivalently, the \emph{out-neighbourhood}, $N^+(x)$); or \label{PropIO}

\vspace{-4.5mm}

\item $N^-(x)$ and $N^+(x) $ separately; or \label{PropIOS}
\item the union $N^-(x) \cup N^+(x)$. \label{PropIOT}
\end{enumerate}


Since our goal is to continue the study of oriented colourings, we
will consider the case where $H$ is irreflexive and the case where
$H$ is reflexive. 
The three definitions are different when $H$ is {reflexive.
When $H$ is irreflexive, definitions \ref{PropIOS} and \ref{PropIOT} coincide. 
Note that, if $H$ is  irreflexive, then $G$ must also be irreflexive.
We  assume throughout the paper that $G$ is an irreflexive oriented graph,
whether or not the target $H$ is reflexive.

In each of the five situations identified above we define a 
\emph{locally-injective oriented $k$-colouring} of an (irreflexive) oriented graph $G$ as a 
locally-injective homomorphism of $G$ to an oriented graph $H$ 
on $k$ vertices.
As noted above, when the locally-injective colouring is proper, $H$ can be taken to be a tournament (irreflexive).
When it is not necessarily proper, $H$ can be taken to be a reflexive tournament.
Thus, studying the complexity of locally-injective oriented colourings leads naturally to
studying the complexity of locally-injective homomorphisms to tournaments and
to reflexive tournaments.



Locally-injective homomorphisms (as in definition \ref{PropIO}) and colourings of oriented graphs were first introduced as an example in monadic second order logic \cite{courcelle}.  Consequently, by Courcelle's Theorem, these problems are all solvable in polynomial time when the input has bounded treewidth.  The same holds for the other definitions above.

Dichotomy theorems are known for each of the five locally-injective 
oriented
$k$-colouring problems \cite{russell,ccm,mrs,mrs1,cobusthesis}.
When the target, $H$, is  reflexive, each such problem 
is NP-complete when $k \geq 3$ and polynomially solvable when $k \leq 2$.
When $H$ is irreflexive, 
each  problem 
is NP-complete when $k \geq 4$ and polynomially solvable when $k \leq 3$.

A fairly complete theory of locally-injective homomorphisms under 
definition \ref{PropIO} has been developed \cite{mrs,mrs2,mrs3,cobusthesis}.   
When the target, $H$,  is reflexive there is a dichotomy theorem characterizing the oriented graphs $H$ for which the problem of deciding the existence of a homomorphism to $H$ is solvable in polynomial time, and those for which it is NP-complete.  When $H$ is irreflexive the complexity has been determined when $H$ has maximum in-degree $\Delta^- \geq 3$ or $\Delta^- \leq 1$; when $\Delta^-=2$ the situation is as rich as that for all digraph homomorphism problems, and hence all constraint satisfaction problems \cite{mrs1}.  
In each case, $H$, for which 
one of these locally-injective homomorphism problems is solvable in polynomial time, it is possible to find a duality theorem for locally-injective homomorphism to $H$.
Duality theorems are also known in cases where the associated locally-injective oriented $k$-colouring problem is polynomial.  These lead to characterization of the critical oriented graphs, i.e., those for which any proper subgraph can be coloured using fewer colours.

Dichotomy theorems for locally-injective homomorphisms to reflexive
tournaments under definitions \ref{PropIOS} and \ref{PropIOT} are presented in \cite{bard}. 
In particular, in each case it is shown that the associated homomorphism problems are NP-complete 
for tournaments on at least three vertices, and solvable in polynomial time for tournaments on at most two vertices. 
With respect to locally-injective homomorphisms to (irreflexive)
tournaments under definition \ref{PropIOS} (and hence also definition \ref{PropIOT}),
it is shown in \cite{ccm}  that the associated homomorphism problems are polynomially 
solvable for tournaments on at most three vertices.

\subsection{Overview}

In this paper we identify obstructions to the polynomially solvable 
locally-injective oriented colouring problems under definitions \ref{PropIOS} and \ref{PropIOT}.
We first identify obstructions to locally-injective homomorphisms to small tournaments and
small reflexive tournaments, and then combine these results.
More notation and terminology is defined in the next section.
Section 3 is concerned with proper locally-injective colourings (i.e.~the target oriented graphs are irreflexive)
under definition \ref{PropIOS} (and hence also definition \ref{PropIOT}).
The only non-trivial case is when $k=3$.
We first describe obstructions as above to
locally-injective homomorphisms to the directed 3-cycle and the transitive tournament on 3 vertices, 
and then combine these to obtain a set of obstructions for locally-injective 3-colouring.
The next two sections, in turn,  describe obstruction sets for improper locally-injective colourings 
(i.e.~the target oriented graphs are reflexive, so adjacent vertices may be assigned the same colour) 
under definitions \ref{PropIOS} and \ref{PropIOT}.
In the latter of these we also give a structural description of the oriented graphs that have a 
locally-injective 2-colouring.





\section{Further notation and terminology}
\label{Sec2}


Throughout the remainder of the paper, the superscript ``$r$'', as in $G^r$, indicates that the oriented graph under consideration is reflexive.  
Oriented graphs without this superscript, as in $G$, are irreflexive.

The \emph{converse} of an oriented graph $G$ is the oriented graph $G^c$ with the same vertex set as $G$, and arc set $\{yx: xy\in E(G)\}$.  

We use $P_n, C_n$, and $T_n$  to denote the directed path on $n$ vertices, the directed cycle on $n$ vertices, and the transitive tournament on $n$ vertices, respectively, $n \geq 1$.  It will be assumed throughout that $C_3$ has vertex set $\{c_1, c_2, c_3\}$ and arc set $\{c_1c_2, c_2c_3, c_3c_1\}$, and that $T_n$ has vertex set $\{t_0, t_1, \ldots, $ $t_{n-1}\}$ and arc set $\{t_it_j: i < j\}$.

We call a homomorphism $f$ of an oriented graph $G$ to an oriented graph $H$:
\begin{itemize}
\item \emph{ios-injective} if, for every vertex $x$ of $G$, the restriction of $f$ to $N^-(x)$ is injective, as is the restriction of $f$ to $N^+(x)$; and
\item   \emph{iot-injective} if, for every vertex $x$ of $G$, the restriction of $f$ to $N^-(x) \cup N^+(x)$ is injective.
\end{itemize}
These two concepts are the same when $H$ is an irreflexive oriented graph, and different when $H$ is a reflexive oriented graph.

The designations ``ios'' and ``iot'' arise from the local injectivity being on \underline{{\bf i}}n-neighbour- hoods and \underline{{\bf o}}ut-neighbourhoods \underline{{\bf s}}eparately, and on \underline{{\bf i}}n-neighbourhoods and \underline{{\bf o}}ut-neigh- bourhoods \underline{{\bf t}}ogether.  In introducing the designations ``ios'' and ``iot'',  the qualifier ``locally'' has been dropped as it is part of the definition. 


It is easy to see that the composition of two ios-injective homomorphisms is an ios-injective homomorphism, and similarly for iot-injective homomorphisms. 


The following structure and its converse will be particularly useful.
We define the \textit{hat} $H_3$ to be the oriented graph with vertex set $V(H_3) = \{v_0, v_1, v_2\}$ 
and edge set $E(H_3) = \{v_0v_1, v_2v_1\}$.   The vertices $v_0$ and $v_2$ will be referred to as the \emph{ends} of $H_3$ or $H_3^c$; $v_1$ is the \emph{apex}.  Whether or not $H$ is reflexive, in an ios-injective or iot-injective homomorphism of $H_3$ or $H_3^c$ to $H$, the vertices $v_0$ and $v_2$ must have different images.

\begin{figure}[htbp]
\begin{center}
\begin{tikzpicture}
	\begin{pgfonlayer}{nodelayer}
		\node [style=vertex] (0) at (-5, 1) {};
		\node [style=vertex] (1) at (-4.5, 2) {};
		\node [style=vertex] (2) at (-4, 1) {};
		\node [style=vertex] (3) at (-2, 1) {};
		\node [style=vertex] (4) at (-1.5, 2) {};
		\node [style=vertex] (5) at (-1, 1) {};
		\node [style=vertex] (6) at (-4.5, 2) {};
		\node [style=box] (7) at (-4.5, 0.25) {$H_3$};
		\node [style=box] (8) at (-1.5, 0.25) {$H_3^c$};
	\end{pgfonlayer}
	\begin{pgfonlayer}{edgelayer}
		\draw [style=arc] (0) to (1);
		\draw [style=arc] (2) to (1);
		\draw [style=arc] (4) to (3);
		\draw [style=arc] (4) to (5);
	\end{pgfonlayer}
\end{tikzpicture}
\caption{The hat and its converse.}
\label{fig1}
\end{center}
\end{figure}
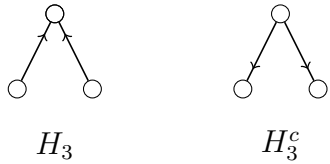

\section{Ios-injective and iot-injective homomorphisms to irreflexive targets}
\label{HomSection}

In this section, we consider ios-injective homomorphisms, and hence also iot-injective homomorphisms, to the four (irreflexive) tournaments on at most three vertices, namely $T_1, T_2, T_3$, and $C_3$.  Since the target oriented graphs have no loops, adjacent vertices of the input graph must be assigned different images. 

Clearly, an oriented graph $G$ has an ios-injective homomorphism to $T_1$ if and only if it has no arcs.  Equivalently, if and only if $T_2$ has no (ios-injective) homomorphism to $G$, if and only if $T_2$ is not a subgraph of $G$.

It is also clear that an oriented graph has an ios-injective homomorphism to $T_2$ if and only if it contains no oriented path on three vertices, that is, if and only if none of $H_3, H_3^c$, and $P_3$ is a subgraph of $G$.  Equivalently, if and only if none of these has an ios-injective homomorphism to $G$.

An oriented graph $G$ with an ios-injective homomorphism to $C_3$ has maximum in-degree at most one and maximum out-degree at most one.  The degree conditions imply that neither $H_3$ nor $H_3^c$ has an ios-injective homomorphism to $G$, so $G$ is a union of directed paths and directed cycles of length a multiple of three.  Any homomorphism of such a graph to $C_3$ is ios-injective. We therefore have the proposition which follows.

\begin{prop} \label{c3} An oriented graph $G$ has an ios-injective homomorphism to $C_3$ if and only if none of $H_3, H_3^{c}$, and $C_n, n \not \equiv 0$ (mod 3) has an injective homomorphism to $G,$ if and only if none of these oriented graphs is a subgraph of $G$.
\end{prop}

Finally, we turn our attention to ios-injective homomorphisms to $T_3$.   The situation is a bit more complicated.   The directed path $P_4$ is an obstruction to ios-injective homomorphism to $T_3$, as are the directed graphs $H_4$, with vertex set $\{h_0, h_1, h_2, h_3\}$ and edge set $\{h_0h_3, h_1h_3, h_2h_3\}$, and $A_4$, with vertex set $\{h_0, h_1, h_2, h_3\}$ and edge set $\{h_0h_2, h_1h_2, h_2h_3\}$. The oriented graphs $H_4^c$  and $A_4^c$ are also
obstructions. Further, we must include as obstructions the graphs obtained by joining two directed paths of length two at the initial vertices, say $H_5$, and its converse, $H_5^c$.  Finally, let $B_2$ be the oriented graph obtained from a path of length two by replacing each edge by a copy of $H_3$. See Figure \ref{fig2}.

\medskip

\begin{figure}[htbp]
\begin{center}
\begin{tikzpicture}
	\begin{pgfonlayer}{nodelayer}
		\node [style=vertex] (0) at (-4, 7) {};
		\node [style=vertex] (1) at (-4, 6) {};
		\node [style=vertex] (2) at (-3, 6) {};
		\node [style=vertex] (3) at (-3, 7) {};
		\node [style=vertex] (4) at (-1, 6) {};
		\node [style=vertex] (5) at (0, 6) {};
		\node [style=vertex] (6) at (-1, 7) {};
		\node [style=vertex] (7) at (0, 7) {};
		\node [style=box] (8) at (-3.5, 5.5) {$H_4$};
		\node [style=box] (9) at (-0.5, 5.5) {$H_4^c$};
		\node [style=vertex] (10) at (-4, 4.5) {};
		\node [style=vertex] (11) at (-3, 4.5) {};
		\node [style=vertex] (12) at (-3.5, 3.75) {};
		\node [style=vertex] (13) at (-3.5, 2.75) {};
		\node [style=vertex] (14) at (-1, 4.5) {};
		\node [style=vertex] (15) at (0, 4.5) {};
		\node [style=vertex] (16) at (-0.5, 3.75) {};
		\node [style=vertex] (17) at (-0.5, 2.75) {};
		\node [style=box] (18) at (-3.5, 2.25) {$A_4$};
		\node [style=box] (19) at (-0.5, 2.25) {$A_4^c$};
		\node [style=vertex] (20) at (-4, 1.25) {};
		\node [style=vertex] (21) at (-4, 0.5) {};
		\node [style=vertex] (22) at (-4, -0.25) {};
		\node [style=vertex] (23) at (-3.25, -0.25) {};
		\node [style=vertex] (24) at (-2.5, -0.25) {};
		\node [style=vertex] (25) at (-1, 1.25) {};
		\node [style=vertex] (26) at (-1, 0.5) {};
		\node [style=vertex] (27) at (-1, -0.25) {};
		\node [style=vertex] (28) at (-0.25, -0.25) {};
		\node [style=vertex] (29) at (0.5, -0.25) {};
		\node [style=box] (30) at (-3.5, -0.75) {$H_5$};
		\node [style=box] (31) at (-0.5, -0.75) {$H_5^c$};
		\node [style=vertex] (32) at (-4.25, -2.5) {};
		\node [style=vertex] (33) at (-3.75, -1.75) {};
		\node [style=vertex] (34) at (-3.25, -2.5) {};
		\node [style=vertex] (35) at (-2.75, -1.75) {};
		\node [style=vertex] (36) at (-2.25, -2.5) {};
		\node [style=vertex] (37) at (-0.25, -2.5) {};
		\node [style=vertex] (38) at (0.25, -1.75) {};
		\node [style=vertex] (39) at (-1.25, -2.5) {};
		\node [style=vertex] (40) at (0.75, -2.5) {};
		\node [style=vertex] (41) at (-0.75, -1.75) {};
		\node [style=box] (42) at (-3.25, -3) {$B_2$};
		\node [style=box] (43) at (-0.25, -3) {$B_2^c$};
	\end{pgfonlayer}
	\begin{pgfonlayer}{edgelayer}
		\draw [style=arc] (0) to (1);
		\draw [style=arc] (3) to (1);
		\draw [style=arc] (2) to (1);
		\draw [style=arc] (4) to (6);
		\draw [style=arc] (4) to (7);
		\draw [style=arc] (4) to (5);
		\draw [style=arc] (10) to (12);
		\draw [style=arc] (11) to (12);
		\draw [style=arc] (12) to (13);
		\draw [style=arc] (17) to (16);
		\draw [style=arc] (16) to (14);
		\draw [style=arc] (16) to (15);
		\draw [style=arc] (22) to (21);
		\draw [style=arc] (21) to (20);
		\draw [style=arc] (22) to (23);
		\draw [style=arc] (23) to (24);
		\draw [style=arc] (25) to (26);
		\draw [style=arc] (26) to (27);
		\draw [style=arc] (29) to (28);
		\draw [style=arc] (28) to (27);
		\draw [style=arc] (32) to (33);
		\draw [style=arc] (34) to (33);
		\draw [style=arc] (34) to (35);
		\draw [style=arc] (36) to (35);
		\draw [style=arc] (41) to (39);
		\draw [style=arc] (41) to (37);
		\draw [style=arc] (38) to (37);
		\draw [style=arc] (38) to (40);
	\end{pgfonlayer}
\end{tikzpicture}
\caption{Some obstructions to injective homomorphism to $T_3$.}
\label{fig2}
\end{center}
\end{figure}
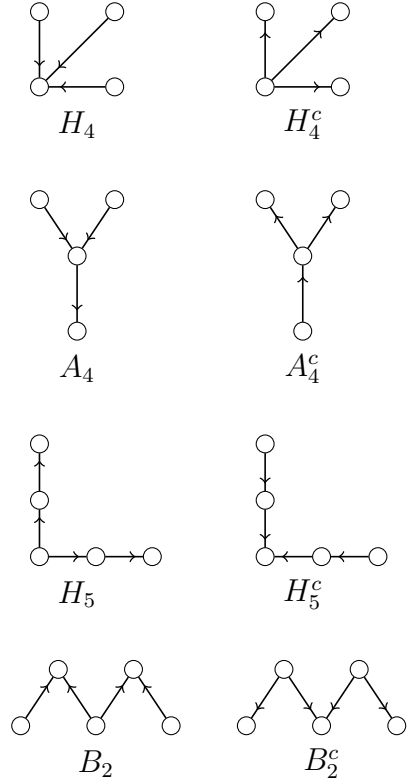

Note that forbidding $H_4, H_4^c, A_4$, and $A_4^c$ means forbidding vertices of degree at least 3 so that the corresponding graphs are simply unions of oriented cycles and oriented paths.

\begin{thm} \label{t3} An oriented graph $G$ has an ios-injective homomorphism to $T_3$ if and only if none of $P_4, H_4, H_4^{c}, H_5, H_5^{c}, A_4, A_4^{c}, B_2 $, and $B_2^{c}$ has an ios-injective homomorphism to $G$.
\end{thm}

\noindent \textbf{Proof.} Clearly none of the given oriented graphs has an ios-injective homomorphism to $T_3$ and hence not to $G$. Suppose none has an ios-injective homomorphism to $G$. We will show that there is an ios-injective homomorphism of $G$ to $T_3$. It can be assumed that $G$ is connected.

Suppose first that $G \not\cong H_3$ has no directed path of length two. (Clearly there is an ios-injective homomorphism of $H_3$ to $T_3$.) Then every vertex has either in-degree zero or out-degree zero. Since neither of $H_4$ and $H_4^c$ has an ios-injective homomorphism to $G$, vertices of in-degree zero have maximum out-degree two and vertices of out-degree zero have maximum in-degree two. Let $S$ be the set of vertices of in-degree zero. Map the vertices in $S$ of out-degree two to $t_0$ and those of out-degree one to $t_1$. Now, construct an undirected graph $L$ with $V(L) = V(G) \setminus S$ and $uv \in E(L)$ when $u$ and $v$ are the end vertices of $H_3^c$ in $G$. Since $B_2^c$ has no ios-injective homomorphism to $G$, the (undirected) graph $L$ contains no path of length two and hence each component of $L$ is isomorphic to $K_1$ or $K_2$.  A 2-colouring of $L$  with the colours $t_1$ and $t_2$ gives an ios-injective homomorphism of $G$ to $T_3$. In particular, let $c$ be the 2-colouring of $L$ such that vertices of $L$ that have in-degree two in $G$ receive color $t_2$ and all remaining vertices, namely those that have in-degree one in $G$, receive colour $t_1$.

Now suppose that $G$ has a directed path of length two. The vertices of any such path must be  coloured $t_0, t_1$ and $t_2$, in order. As before, let $S$ be the set of vertices of in-degree zero. Because none of the given oriented graphs has an ios-injective homomorphism to $G$, the set $S$ contains all vertices of out-degree two, all vertices that are the initial vertex of a directed path of length two, and all vertices of out-degree one (and in-degree zero) that are an endpoint of a hat $H_3$. Further, $S$ is an independent set.  Every vertex not in $S$ is the apex of a hat $H_3$, the end of a hat $H_3^c$, or a midpoint or terminus of a directed path of length two. Again, construct an undirected graph $L$ with $V(L) = V(G) \setminus S$ and $uv \in E$ when $u$ and $v$ are the end vertices of $H_3^c$ in $G$. As before, each component of $L$ is isomorphic to $K_1$ or $K_2$.

Let $X$ be the set of terminal vertices of directed paths of length two in $G$ together with the set of apices of hats $H_3$, and let $Y$ be the set of midpoints of directed paths of length two. 
Having neither of $A_4$ and $A_4^c$ as a preimage implies that $Y$ is independent and having no $H_5^c$ as a preimage implies that $X$ is independent.
Further, having no $P_4$ as a preimage implies that the sets $X$ and $Y$ are disjoint. Given the simple structure of $L$, it has a 2-colouring $c$ in which the vertices in $X$ are assigned colour $t_2$ and the vertices in $Y$ are assigned colour $t_1$. We now need only assign colours to those vertices that are one end of an $H_3^c$. If the other end is also the midpoint of a path of length two, assign colour $t_2$. Otherwise, the other end is the apex of an $H_3$; assign colour $t_1$.

In this case, an ios-injective homomorphism of $G$ to $T_3$ is obtained by mapping the vertices in $S$ of out-degree one that are ends of an $H_3$ to $t_1$ and all other vertices in $S$ to $t_0$, and each vertex of $L$ to the vertex of $T_3$ corresponding to its colour. \hfill $\square$

\bigskip

The above theorem leads to a forbidden subgraph characterization of the oriented graphs that have an ios-injective homomorphism to $T_3$.  

\begin{thm}
An oriented graph $G$ has an ios-injective homomorphism to $T_3$ if and only if $G$ has no oriented 4-cycle and none of $C_3, P_4, H_4, H_4^c, H_5, H_5^c, A_4$, $A_4^c, B_2$ and $B_2^c$ is a subgraph of $G$. \label{iosT_3forb} 
\end{thm}

\noindent \textbf{Proof.} 
Any oriented graph among $C_4, H_4, H_4^c, A_4$ and $A_4^c$ has an ios-injective homomorphism to an oriented graph $G$ if and only if it is a subgraph of $G$.  If one of $H_5, H_5^c, B_2, B_2^c$ and $P_4$ has an ios-injective homomorphism to $G$, then either it is a subgraph of $G$, or some other oriented graph in the collection (including the 4-cycle) is a subgraph of $G$. \hfill $\square$

\bigskip
Since the list of forbidden subgraphs is finite, it follows that there is a polynomial time algorithm to decide if a given oriented graph $G$ has an ios-injective homomorphism to $T_3$.  In cases where there is no such homomorphism, a forbidden subgraph is found.  The proof of Theorem \ref{t3} gives an efficient method for finding the ios-injective homomorphism when it exists.

\section{Ios-injective homomorphisms to reflexive targets}

In this section, we consider ios-injective homomorphisms to the two reflexive tournaments on at most two vertices.  

It is clear that an oriented graph $G$ has an ios-injective homomorphism to $T_1^r$ if and only if it has maximum in-degree at most one and maximum out-degree at most one, that is, if and only if neither $H_3$ nor $H_3^c$ has an ios-injective homomorphism to $G$.  The last condition is equivalent to neither of these oriented graphs being a subgraph of $G$.  Consequently, the only oriented graphs which have an ios-injective homomorphism to $T_1^r$ are disjoint unions of directed paths and directed cycles.

\subsection{Ios-injective homomorphism to $T_2^r$}
\label{Sec4}

%
%
In this sub-section we show that the existence of an ios-injective homomorphism to $T_2^r$ is solvable in polynomial time by
describing a finite list of substructures which preclude the existence of such a mapping.  
The approach is to set up a colouring extension problem, and extract the obstructions from situations where no extension exists.  This method has been used elsewhere (e.g.~see \cite{mrs2,mrs1}).

Let $X_2$ denote the orientation of the complete bipartite graph $K_{1, 4}$ with a vertex having in-degree two and out-degree two. We clearly have:

\begin{prop}
None of the oriented graphs $H_4, H_4^c,
X_2$ has an ios-injective homomorphism to $T_2^r$. 
\label{ios-forb}
\end{prop}

Let $G$ be an oriented graph that has maximum in-degree two and maximum out-degree two (note: $H_4$ and $H_4^c$ are among the forbidden images), and such that none of the oriented graphs in Proposition \ref{ios-forb} has an ios-injective homomorphism to $G$.  

Let $\mathcal{O}_2$ be the set of vertices with out-degree 2, and let $\mathcal{I}_2$
be the set of vertices with in-degree 2.
Since $X_2$ does not have an ios-injective homomorphism to $G$, 
$\mathcal{I}_2 \cap \mathcal{O}_2 = \emptyset$.

Every vertex in $\mathcal{I}_2$ must map to $t_1$, 
and every vertex in $\mathcal{O}_2$ must map to $t_0$.  
The vertices in $V - (\mathcal{I}_2 \cup \mathcal{O}_2)$
all have in-degree at most 1 and out-degree at most 1.
These must lie on directed paths which are either components of $G$, or start or end at 
vertices in $\mathcal{I}_2 \cup \mathcal{O}_2$.

If a vertex $x$ maps to $t_1$, then any vertex reachable from $x$ by a directed path must map to $t_1$.  
Similarly, if $y$ maps to $t_0$, then any vertex that can reach $y$ by a directed path must also map to $t_0$. 
Finally, vertices that are ends of a copy of $H_3$ or $H_3^c$ cannot have the same image. 
These ideas lead to the version of the {\em Consistency Check Algorithm} given below. 
 In the algorithm, the list $\ell(x)$ corresponding to a vertex $x$ contains its possible images, and a possible image that cannot be assigned to $x$ gets removed from $\ell(x)$.
 
\vskip 0.25 in
\noindent{\bf Consistency Check Algorithm}\\ \noindent
Let $G$ be an oriented graph that has maximum in-degree two and maximum out-degree two such that none of the oriented graphs in Proposition \ref{ios-forb} has an ios-injective homomorphism to $G$.\\ \\
Initialization\\ 
$$\ell(x) = \begin{cases}
\{t_0\} & \mathrm{if}\ d^+(x) = 2;\\
\{t_1\} & \mathrm{if}\ d^-(x) = 2;\\
\{t_0, t_1\} & \mathrm{otherwise}.\\
\end{cases} $$
\noindent

\noindent Repeat until some list becomes empty or no lists change.
\begin{enumerate}
\item If $\ell(y) = \{t_0\}$ and $xy \in E$,  then remove $t_1$ from $\ell(x)$;
\item  If $\ell(y) = \{t_1\}$ and $yx \in E$,  then remove $t_0$ from $\ell(x)$;
\item If $\ell(y) = \{c\}$ and $xs, ys \in E$ or $sx, sy \in E$, for some $s \in V(G)$, then remove $c$ from $\ell(x)$.
\end{enumerate}

\bigskip
If the Consistency Check Algorithm terminates with no list empty, we say it \emph{succeeds}.  
Otherwise,  we say  it \emph{fails}.
If the Consistency Check Algorithm fails then $G$ has no ios-injective homomorphism to $T_2^r$ because
there is a vertex that cannot be assigned an image.  

\begin{prop}
Let $G$ be an oriented graph that has maximum in-degree two and maximum out-degree two  such that none of the oriented graphs in Proposition \ref{ios-forb} has an ios-injective homomorphism to $G$.  
The Consistency Check Algorithm succeeds when run on $G$ if and only if $G$ has an ios-injective homomorphism to 
$T_2^r$.
\label{CCsucceeds}
\end{prop}
\noindent\textbf{Proof}.
If the Consistency Check Algorithm fails then $G$ has no ios-injective homomorphism to $T_2^r$ because
there is a vertex that cannot be assigned an image.  Suppose, then, that the algorithm succeeds.

Suppose first that every list contains exactly one element when the algorithm terminates.
Then, by the instructions of the algorithm, mapping each vertex of $G$ to the element in its list
is an ios-injective homomorphism to $T_2^r$.

Suppose the lists of some vertices contain two elements when the algorithm terminates.
Then these vertices lie on directed paths which are either components of $G$, or start at vertices 
of $\mathcal{O}_2$ or end at vertices of $\mathcal{I}_2$.

Let $P$ be a directed path which is itself a component of $G$.
Since no vertex of $P$ belongs to $\mathcal{I}_2 \cup \mathcal{O}_2$, 
 the lists corresponding to vertices of $P$ are not changed
by the Consistency Check Algorithm.
Every vertex of $P$ can be mapped to $t_1$. 

Let $Q$ be a directed path that starts at a vertex of $\mathcal{O}_2$.
Then $Q$ either ends at a vertex of $\mathcal{I}_2$ or at a vertex of out-degree 0.
In the latter case, all vertices of $Q$ can be mapped to $t_0$.
The case of a directed path that starts at a vertex of in-degree 0 and ends at a vertex of $\mathcal{I}_2$ is similar.

It remains to assign images to the vertices on directed paths from vertices in 
$\mathcal{O}_2$ to vertices in $\mathcal{I}_2$.
Recall that $\mathcal{O}_2 \cap \mathcal{I}_2 = \emptyset$.
Construct a bipartite graph $B$ with bipartition $(\mathcal{O}_2, \mathcal{I}_2)$
such that $xy \in E(B)$ if and only if $x \in \mathcal{O}_2$, $y \in \mathcal{I}_2$ and
there is a directed path in $G$ from $x$ to $y$ in which every internal vertex has
a list of size 2.
Since $G$ contains neither $H_4$ nor $H_4^c$, the graph $B$ has maximum 
degree 2.
Therefore, $B$ has a 2-edge colouring using the colours $t_0$ and $t_1$.
Assign the colour of each edge $uw \in E(B)$ to all internal vertices 
of the corresponding directed path in $G$.
This gives an  ios-injective 
homomorphism of $G$ to $T_2^r$.



Hence, if the Consistency Check Algorithm succeeds, then $G$ has an ios-injective 
homomorphism to $T_2^r$.
$\Box$

In order to identify the remaining obstructions we consider situations when the Consistency Check Algorithm fails.
%

In what follows, for $i, j \in \{0, 1\}$, the term {\em $(i,j)$-walk from $x$ to $y$} indicates a sequence of adjacent vertices such that if $x$ maps to $t_i$ then $y$ must map to $t_j$.  The walks defined below consist of oppositely oriented directed paths $F_i$ and $B_j$.  The arcs of each $F_i$ are oriented in the \emph{forwards} direction, that is, the tail of each arc precedes its head in the walk.  The arcs of each $B_j$ are oriented in the \emph{backwards} direction, meaning that the head of each arc precedes its tail in the walk (see Figure~\ref{FigWalks}).   

A {\em (1,0)-walk from $x$ to $y$} is either
the subgraph $H_3^c$ with vertices $x, s,$ and $y$, or a sequence 
$$F_1, B_1, F_2, B_2, \ldots, F_{k}, B_{k},$$
where: 
\begin{enumerate}
\item $k \geq 1$;
\item $x$ is the initial vertex of $F_1$, and $y$ is the initial vertex of $B_{k}$;
\item the forwards directed path $F_1$ has length at least 1;
\item each of the directed paths  $B_1, F_2, B_2, \ldots, F_{k}$ has length at least 2;
\item the final backwards directed path $B_k$ has length at least 1.
\end{enumerate}

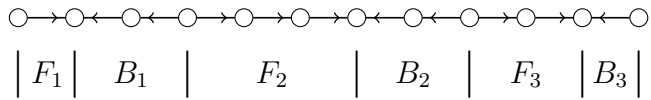
\begin{figure}[htbp]
\begin{center}
\begin{tikzpicture}
	\begin{pgfonlayer}{nodelayer}
		\node [style=vertex] (0) at (-6, 0) {};
		\node [style=vertex] (1) at (-5.25, 0) {};
		\node [style=vertex] (2) at (-4.5, 0) {};
		\node [style=vertex] (3) at (-3.75, 0) {};
		\node [style=vertex] (4) at (-3, 0) {};
		\node [style=vertex] (5) at (-2.25, 0) {};
		\node [style=vertex] (6) at (-1.5, 0) {};
		\node [style=vertex] (7) at (-0.75, 0) {};
		\node [style=vertex] (8) at (0, 0) {};
		\node [style=vertex] (9) at (0.75, 0) {};
		\node [style=vertex] (10) at (1.5, 0) {};
		\node [style=vertex] (11) at (2.25, 0) {};
		\node [style=box, scale=1.5] (12) at (-6, -0.75) {$|$};
		\node [style=box, scale=1.5] (13) at (-5.25, -0.75) {$|$};
		\node [style=box, scale=1.5] (14) at (-3.75, -0.75) {$|$};
		\node [style=box, scale=1.5] (15) at (-1.5, -0.75) {$|$};
		\node [style=box, scale=1.5] (16) at (0, -0.75) {$|$};
		\node [style=box, scale=1.5] (17) at (1.5, -0.75) {$|$};
		\node [style=box, scale=1.5] (18) at (2.25, -0.75) {$|$};
		\node [style=box] (19) at (-5.61, -0.75) {$F_1$};
		\node [style=box] (20) at (-4.5, -0.75) {$B_1$};
		\node [style=box] (21) at (-2.625, -0.75) {$F_2$};
		\node [style=box] (22) at (-0.75, -0.75) {$B_2$};
		\node [style=box] (23) at (0.75, -0.75) {$F_3$};
		\node [style=box] (24) at (1.87, -0.75) {$B_3$};
	\end{pgfonlayer}
	\begin{pgfonlayer}{edgelayer}
		\draw [style=arc] (0) to (1);
		\draw [style=arc] (2) to (1);
		\draw [style=arc] (3) to (2);
		\draw [style=arc] (3) to (4);
		\draw [style=arc] (4) to (5);
		\draw [style=arc] (5) to (6);
		\draw [style=arc] (8) to (7);
		\draw [style=arc] (7) to (6);
		\draw [style=arc] (8) to (9);
		\draw [style=arc] (11) to (10);
		\draw [style=arc] (9) to (10);
	\end{pgfonlayer}
\end{tikzpicture}
\caption{A $(1,0)$-walk $F_1, B_1, F_2, B_2, F_3, B_3$.}
\label{FigWalks}
\end{center}
\end{figure}

Similarly, we define $(0,1)$-walks, $(0,0)$-walks, and $(1,1)$-walks.  In what follows, $F_i$ denotes a forwards directed path, $B_j$ denotes a backwards directed path, $k \geq 1$,  each directed path $F_i$ or $B_j$ other then the first and last has length at least 2, and the first and last directed paths have length at least 1.
\begin{itemize}
\item A {\em $(0,1)$-walk} from $x$ to $y$ is either the subgraph $H_3$ with vertices $x, z,$ and $y$ or a sequence
$B_1, F_1, F_2, B_2, \ldots, B_{k}, F_k$,  where $x$ is the terminal vertex of $B_1$ and $y$ is the terminal vertex of $F_{k}$.

\item A {\em $(0,0)$-walk} from $x$ to $y$ is a sequence $B_1, F_1, F_2, B_2, \ldots, B_{k}$, where  $x$ is the terminal vertex of $B_1$ and $y$ is the initial vertex of $B_{k}$.

\item A {\em $(1,1)$-walk} from $x$ to $y$ is a sequence $F_1, B_1, F_2, B_2, \ldots, F_{k}$, where  $x$ is the initial vertex of $F_1$ and $y$ is the terminal vertex of $F_{k}$.

\end{itemize}
In all cases, we refer to $x$ as the \emph{first vertex of the walk} and $y$ as the \emph{last vertex of the walk}.

\begin{prop}
Let $G$ be an orientation of a graph of maximum degree 3.
Suppose the Consistency Check Algorithm is run on $G$.
\begin{enumerate}
\item Suppose $\ell(x) = \{t_0\}$ at some point during the Consistency Check Algorithm.
\begin{enumerate}
\item If there exists a $(0,0)$-walk from $x$ to $y$ then $t_1 \not\in \ell(y)$ at the conclusion of the Algorithm;
\item If there exists a $(0,1)$-walk from $x$ to $y$ then $t_0 \not\in \ell(y)$ at the conclusion of the Algorithm.
\end{enumerate}
\item Suppose $\ell(x) = \{t_1\}$ at some point during the Consistency Check Algorithm.
\begin{enumerate}
\item If there exists a $(1,1)$-walk from $x$ to $y$ then $t_0 \not\in \ell(y)$ at the conclusion of the Algorithm;
\item If there exists a $(1,0)$-walk from $x$ to $y$ then $t_1 \not\in \ell(y)$ at the conclusion of the Algorithm.
\end{enumerate}
\end{enumerate}
\label{RemovalProp}
\end{prop}

\noindent \textbf{Proof}.
We prove only statement (2) (b).  The proofs of the other statements are similar.

Suppose $\ell(x) = \{t_1\}$ at some point in the Consistency Check Algorithm.
If the (1, 0)-walk is $H_3^c$ with vertices $x, s, y$ then $t_1$ is removed from $\ell(y)$ when step (3) is applied at $x$.
Suppose, then, that the $(1,0)$-walk is a sequence 
$F_1, B_1, F_2, B_2, \ldots, F_{k}, B_{k}$.

Let $x'$ be the terminal vertex of $F_1$ and $y'$ the in-neighbour of $x'$ that belongs to $B_1$.
Then, by the instructions of the algorithm, $t_0$ is removed from the list of each vertex joined to $x$ by an arc of $F_1$.  
Thus, eventually, $\ell(x') = \{t_1\}$.
Since $F_1$ has length at least 1, there is a copy of $H_3$ consisting of $x'$ and its two in-neighbours.
Thus $t_1$ is removed from $\ell(y')$  and hence eventually from the list of each vertex of $B_1$ joined to
$y'$ by a directed path.
Continuing in this way, $t_0$ is removed from the list of the terminal vertex of 
each forwards directed path, and $t_1$ is removed from the list of the terminal vertex of each 
backwards directed path.
In particular, $t_1$ is removed from $\ell(y)$.
$\Box$

It is natural to wonder about the restriction that only the first and last directed paths can have length 1.
If this is not so, then our walk contains 3 consecutive oppositely oriented arcs, say $ab, cb,$ and $cd$.
Following the argument above, $t_0 \not\in \ell(a)$, so if our oriented graph has an ios-injective homomorphism
to $T_2^r$, then $a$ maps to $t_1$.  Hence $b$ must also map to $t_1$, which is why $t_0$ is removed from $\ell(b)$.
Subsequently, $c$ maps to $t_0$.  
But then both $t_0$ and $t_1$ are legitimate possible images for $d$, so neither can be removed from its list.
The case where the arcs are $ba, bc$, and $dc$ is similar.

Suppose the Consistency Check Algorithm fails when run on $G$.
The arguments toward empty lists necessarily constructs substructures for the obstructions to ios-injective homomorphism to $T_2^r$.
The instructions of the algorithm imply the following:

\begin{prop}\label{eliminate}
Let $G$ be an orientation of a graph with maximum degree 3.
Suppose the Consistency Check Algorithm is run on $G$.
\begin{enumerate}
\item If $t_0$ is deleted from the list of $y$ then one of the following holds:
\begin{enumerate}
\item there is a vertex $x$ with list $\{t_0\}$ and a $(0, 1)$-walk to $y$; or
\item there is a vertex $x$ with list $\{t_1\}$ and a $(1, 1)$-walk to $y$.
\end{enumerate}
\item  If $t_1$ is deleted from the list of $y$ then one of the following holds:
\begin{enumerate}
\item there is a vertex $x$ with list $\{t_0\}$ and a $(0, 0)$-walk to $y$; or
\item there is a vertex $x$ with list $\{t_1\}$ and  a $(1,0)$-walk to $y$.
\end{enumerate}
\end{enumerate}
\end{prop}

For $j \in \{0,1\}$ let $\mathcal{F}_{0,j}$ be the set of all oriented paths 
constructed from a $(0, j)$-walk as follows.
If the walk is neither $H_3$ nor $H_3^c$, 
then add two arcs (and the corresponding different vertices) from
the initial (first) vertex of the walk so it has out-degree 2 and in-degree 1.
If $j = 1$ and the $(0, 1)$-walk is $H_3$, then add an arc and corresponding vertex from its inital vertex so it has out-degree 2.
If $j = 1$ and the $(0, 1)$-walk is $H_3^c$, then add an arc and corresponding vertex to its terminal (last) vertex so it has in-degree 2.
Similarly, for $j \in \{0,1\}$ let $\mathcal{F}_{1,j}$ be the set of all oriented paths 
constructed from a $(1, j)$-walk which is neither $H_3$ nor $H_3^c$ by adding two arcs (and the corresponding different vertices) to the initial vertex of the walk so it has in-degree 2 and out-degree 1.
If $j = 1$ and the $(1, 0)$-walk is $H_3$, then add an arc and corresponding vertex from its terminal vertex so it has out-degree 2.
If $j = 1$ and the $(1, 0)$-walk is $H_3^c$, then add an arc and corresponding vertex to its initial vertex so it has in-degree 2.

\begin{prop}
Let $G$ be an orientation of a graph with maximum degree 3.
For $i, j \in \{0, 1\}$, there is an $(i,j)$-walk in $G$ from $x \in \mathcal{O}_2 \cup \mathcal{I}_2$ to $y$ if and only if there exists $F \in \mathcal{F}_{i,j}$ for which there is an ios-injective homomorphism of $F$ to $G$ such that the first vertex of $F$ maps to $x$ and the last vertex of $F$ maps to $y$.
\end{prop}


Let $G$ be a given oriented graph for which we want to test the existence of an ios-injective homomorphism to $T_2^r$. 
If $\mathcal{O}_2 \cap \mathcal{I}_2 \not= \emptyset$, then $G$ contains $X_2$.   
If $G$ has a vertex of in-degree 3, then $G$ contains $H_4$, and if $G$ has a vertex
of out-degree 3, then $G$ contains $H_4^c$.
 Otherwise, $G$ is an orientation of a graph with maximum degree 3.
Run the Consistency Check Algorithm on $G$.  
By Proposition \ref{CCsucceeds}, the algorithm succeeds if and only if there is an ios-injective homomorphism 
to $T_2^r$.

Suppose the Consistency Check Algorithm fails.  
Then the list of some vertex $y$ becomes empty.
By Proposition \ref{eliminate}, this occurs because for $i_1, i_2 \in \{0, 1\}$
there is both an $(i_1, 0)$-walk and an $(i_2, 1)$-walk for which $y$ is the last vertex.

The lemma which follows is a consequence of Proposition \ref{eliminate}.


\begin{lem}
Let $G$ be an orientation of a graph with maximum degree 3.
If the Consistency Check Algorithm fails when run on $G$, 
then for $p, q \in \{0, 1\}$
there is both a $(p,0)$-walk and a $(q,1)$-walk in $G$ for which $y$ is the last vertex.
\label{nrg}
\end{lem}

Let $\mathcal{F}$ be the set of all oriented graphs constructed as follows.
For  $p, q \in \{0, 1\}$, let $F_0 \in \mathcal{F}_{p, 0}$ and $F_1 \in \mathcal{F}_{q, 1}$.
Identify the last vertex of $F_0$ and the last vertex of $F_1$.

\medskip

The following theorem summarizes our work in this sub-section.

\begin{thm}
An oriented graph $G$ admits an ios-injective 2-colouring if and only if no oriented graph in $\{H_4, H_4^c, X_2\} \cup  \mathcal{F}$ has an ios-injective homomorphism to $G$.
\label{ios-refl-2}
\end{thm}

We conclude the section with an example of finding the obstruction that certifies $T_3$ has no ios-injective homomorphism to $T_2^r$.
Referring to  Figure \ref{iosExample}, when the Consistency Check Algorithm is initialized to run on $T_3$ we have
$\ell(a) = \{t_0\},\ \ell(b) = \{t_0, t_1\}$ and $\ell(c) = \{t_1\}$.
The $(0, 1)$-walk with vertex sequence $a, c, b$ leads to $t_0$ being removed from $\ell(b)$, and corresponds to
an ios-injective homomorphism of  the element of $\mathcal{F}_{0, 1}$ shown on the middle left to $T_3$.
The $(1, 0)$-walk with vertex sequence $c, a, b$ leads to $t_1$ being removed from $\ell(b)$, and corresponds to
an ios-injective homomorphism of the element of $\mathcal{F}_{1, 0}$ shown on the middle right to $T_3$.
Identifying the solid vertex of these two walks  leads to the oriented path $F \in \mathcal{F}$ shown on the right.
It is easy to check that $F$ has no ios-injective homomorphism to $T_2^r$: when the Consistency Check Algorithm is run on $F$, the list of the solid vertex becomes empty.
The labels on the vertices of $F$ describe an ios-injective homomorphism of $F$ to $T_3$.

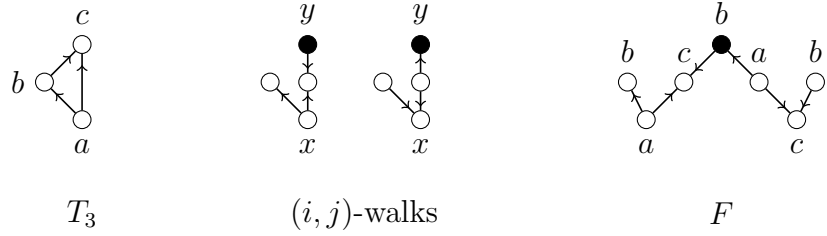
\begin{figure}[htbp]
\begin{center}
\begin{tikzpicture}
	\begin{pgfonlayer}{nodelayer}
		\node [style=vertex, label=below:$a$] (0) at (-6, -0.5) {};
		\node [style=vertex, label=left:$b$] (1) at (-6.5, 0) {};
		\node [style=vertex, label=above:$c$] (2) at (-6, 0.5) {};
		
		\node [style=box] (100) at (-6, -1.75) {$T_3$};
		
		\node[style=vertex, label=below:$x$] (10) at (-3, -0.5) {};
		\node[style=vertex] (11) at (-3, 0) {};
		\node[style=blackvertex, label=above:$y$] (12) at (-3, 0.5) {};
		\node[style=vertex] (13) at (-3.5, 0) {};
		
		\node[style=vertex, label=below:$x$] (20) at (-1.5, -0.5) {};
		\node[style=vertex] (21) at (-1.5, 0) {};
		\node[style=blackvertex, label=above:$y$] (22) at (-1.5, 0.5) {};
		\node[style=vertex] (23) at (-2, 0) {};
		
		\node[style=box] (24) at (-2.25, -1.75) {$(i, j)$\mbox{-walks}};
		
		\node[style=vertex, label=below:$a$] (30) at (1.5, -0.5) {};
		\node[style=vertex, label=above:$b$] (31) at (1.25, 0) {};
		\node[style=vertex, label=above:$c$] (32) at (2, 0) {};
		\node[style=blackvertex, label=above:$b$] (33) at (2.5, 0.5) {};
		\node[style=vertex, label=above:$a$] (34) at (3, 0) {};
		\node[style=vertex, label=above:$b$] (35) at (3.75, 0) {};
		\node[style=vertex, label=below:$c$] (36) at (3.5, -0.5) {};
		
		\node[style=box] (37) at (2.5, -1.75) {$F$};

	\end{pgfonlayer} 
	\begin{pgfonlayer}{edgelayer}
		\draw [style=arc] (0) to (1);
		\draw [style=arc] (1) to (2);
		\draw [style=arc] (0) to (2);
		
		\draw [style=arc] (10) to (11);
		\draw [style=arc] (12) to (11);
		\draw [style=arc] (10) to (13);
		
		\draw [style=arc] (21) to (20);
		\draw [style=arc] (21) to (22);
		\draw [style=arc] (23) to (20);
		
		\draw [style=arc] (30) to (31);
		\draw [style=arc] (30) to (32);
		\draw [style=arc] (33) to (32);
		\draw [style=arc] (34) to (33);
		\draw [style=arc] (34) to (36);
		\draw [style=arc] (35) to (36);

	\end{pgfonlayer}
\end{tikzpicture}
\end{center}
\caption{An example of finding an obstruction.}
\label{iosExample}
\end{figure}

\section{Iot-injective homomorphisms to reflexive targets}

In this section we consider iot-injective homomorphisms to the two reflexive tournaments on at most two vertices. 

It is clear that an oriented graph has an iot-injective homomorphism to $T_1^r$ if and only if it contains no oriented path on three vertices, that is, if and only if it is a disjoint union of copies of $T_1$ and $T_2$.   Curiously, these are the same oriented graphs that have an ios-injective homomorphism to $T_2$.  Put differently, an oriented graph $G$ has an iot-injective homomorphism to $T_1^r$ if and only if none of $H_3, H_3^c, P_3$ have an iot-injective homomorphism to $G$. The last condition is equivalent to  none of these oriented graphs being a subgraph of $G$.  
  
We now turn our attention to $T_2^r$.   No orientation of a graph with a vertex of degree three has an iot-injective homomorphism to $T_2^r$.  Hence  any orientation of $K_{1, 3}$ is an obstruction; call the set of these $\mathcal{S}$.  

The set $\mathcal{P}$ of oriented paths constructed below will be shown to be obstructions.   The following definition is useful.  An \emph{alternating matching} in an oriented path or oriented cycle is a matching that saturates every vertex of degree two in the underlying graph, and whose edges alternate in orientation with respect to the vertex ordering of the path or cycle.  Each path $P \in \mathcal{P}$ has:
\begin{itemize}
\item   $V(P) =\{v_0, v_1, v_2, \ldots, v_{2k}\}$, where $k \geq 2$;
\item  The vertices in $\{v_1, v_2, \ldots, v_{2k-1}\}$ induce a subgraph consisting of two disjoint alternating matchings;
\item The arc joining $v_0$ and $v_1$ has the same orientation (with respect to the beginning and end of the path) as the arc joining $v_2$ and $v_3$; and
\item The arc joining $v_{2k-1}$ and $v_{2k}$ has the same orientation as the arc joining $v_{2k-3}$ and $v_{2k-2}$.
\end{itemize}

\begin{thm}
Let $G$ be an oriented graph.  The following statements are equivalent.
\begin{enumerate}
\item[1.] $G$ has an iot-injective homomorphism to $T_2^r$.
\item[2.] No $F \in \mathcal{S} \cup \mathcal{P}$ has an iot-injective homomorphism to $G$.
\item[3.]  Each component $D$ of $G$ satisfies
\begin{enumerate}
\item[(a)] $D$ is an oriented path on at most four vertices; or
\item[(b)] $D$ is an oriented path with an alternating matching; or
\item[(c)] $D$ is an oriented cycle on $4k$ vertices, where $k \geq 1$, with an alternating matching.
\end{enumerate}
\end{enumerate}
\label{iot-refl-2}
\end{thm}

\noindent \textbf{Proof.}
$(1) \Rightarrow (2)$.  No $F \in \mathcal{S} \cup \mathcal{P}$ has an iot-injective homomorphism to $T_2^r$.  Since iot-injective homomorphisms compose, no such $F$ can have an iot-injective homomorphism to $G$.

$(2) \Rightarrow (3)$.
We prove the contrapositive.  Suppose that $G$ satisfies none of the conditions (a), (b) and (c).  If the underlying undirected graph of $G$ has a vertex of degree at least three, then some element of $\mathcal{S}$ has an iot-injective homomorphism to $G$.  Suppose, then, that the underlying graph of $G$ has maximum degree at most two, and hence each component of $G$ is an oriented path of length at least five, or an oriented cycle.

Let $D$ be a component of $G$ which is an oriented path of length at least five.  Let $M_1$ be the maximum matching in $D$ that contains the first arc of $D$, and let $M_2$ be the matching $E(D) - M_1$.  By hypothesis, neither $M_1$ nor $M_2$ is an alternating matching.  Hence each of these matchings has a pair of arcs at distance one with the same orientation;  say $e_{11}, e_{12} \in M_1$ and $e_{21}, e_{22} \in M_2$, where the section from $e_{11}$ to $e_{22}$ contains $e_{12}$ and $e_{21}$. Choose these pairs so that the size of the section, $P$, from $e_{11}$ to $e_{22}$ is minimized.   By choice of the pairs, the oriented path obtained from $P$ by deleting the first and last vertices consists of two alternating matchings.  Hence $P \in \mathcal{P}$ has an iot-injective homomorphism to $D$. 

Now let $D$ be a component of $G$ which is an oriented even cycle.  Then $D$ consists of two disjoint maximum matchings, neither of which is alternating.  As before, there exists $P \in \mathcal{P}$ which has an iot-injective homomorphism to $D$. 

Finally, let $D$ be a component of $G$ that is an oriented odd cycle.  By parity, there exists a pair of arcs $e_{1}, e_{3}$ at distance one with the same orientation.  Let the arcs of $D$ be, in cyclic order,  $e_1, e_2, \ldots, e_{2k+1}$.  Then the path $Q$ consisting of $2k+1$ consecutive arcs oriented the same as these, followed by three arcs oriented the same as $e_1, e_2$ and $e_3$, has an even number of edges and consists of two disjoint maximum matchings, neither of which is alternating. Clearly $Q$ is in $\mathcal{P}$, and has an iot-injective homomorphism to $D$. 

$(3) \Rightarrow (1)$  It is easy to check that each such graph has an iot-injective homomorphism to $T_2^r$.
\hfill $\square$

\section{Acknowledgement}

The authors are grateful for the contribution of the anonymous referee whose careful reading of the manuscript and helpful comments greatly improved the paper.


\bigskip\bigskip\bigskip
R.J.~Campbell\\
Department of Computing Science, University of the Fraser Valley, Abbotsford , BC\\

N.E.~Clarke$^\dag$\\ 
Department of Mathematics and Statistics,  Acadia University, Wolfville, NS\\

G. MacGillivray$^\dag$\\
Department of Mathematics and Statistics, University of Victoria, Victoria, BC\\

\vfill
\hrule width 2.5in
$ ^\dag$ Research supported by NSERC

\end{document}